\documentclass[12pt]{article}
\usepackage{amssymb,amsthm,amscd,array}
\usepackage[final]{epsfig}
\usepackage{graphics}
\usepackage{amsmath}
\usepackage{amsfonts}
\usepackage{latexsym}
\usepackage{amssymb}
\usepackage{graphicx}
\usepackage{epstopdf}
\let\ssection=\section
\renewcommand{\section}{\setcounter{equation}{0}\ssection}

\newtheorem{lemma}{Lemma}[section]
\newtheorem{proposition}[lemma]{Proposition}
\newtheorem{remark}[lemma]{Remark}
\newtheorem{example}[lemma]{Example}
\newtheorem{theorem}{Theorem}
\newtheorem{fact}[lemma]{Fact}
\newtheorem{definition}[lemma]{Definition}
\newtheorem{corollary}[lemma]{Corollary}

\begin{document}

\newcommand{\eps}{{\varepsilon}}
\newcommand{\g}{{\gamma}}
\newcommand{\de}{{\delta}}
\def\proof{\paragraph{Proof.}}
\newcommand{\proofend}{$\Box$\bigskip}
\newcommand{\const}{{\mathrm{const}}}
\newcommand{\C}{{\mathbb C}}
\newcommand{\CP}{{\mathbb CP}}
\newcommand{\cE}{{\mathcal E}}
\newcommand{\cH}{{\mathcal H}}
\newcommand{\PGL}{{\mathrm{PGL}}}
\newcommand{\Sl}{{\mathrm{sl}}}
\newcommand{\Q}{{\mathbb Q}}
\newcommand{\R}{{\mathbb R}}
\newcommand{\T}{{\mathbb T}}
\newcommand{\Z}{{\mathbb Z}}
\newcommand{\RP}{{\mathbb {RP}}}

\def\a{\alpha}
\def\b{\beta}
\def\l{\lambda}
\def\m{\mu}

\title{Hyperbolic Carath\'eodory conjecture}

\author{Valentin Ovsienko and Serge Tabachnikov\\ }
\date{}
\maketitle

\hfill{Dedicated to Vladimir Igorevich Arnold} 

\hfill{in occasion of his 70th anniversary}

\vskip 1cm

\begin{abstract}
A quadratic point on a surface in $\RP^3$
is a point at which the surface can
be approximated by a quadric abnormally well (up to order 3).
We  conjecture that the least number of quadratic points 
on a generic compact non-degenerate hyperbolic surface is 8;
the relation between this and the classic
Carath\'eodory conjecture is similar to the relation between 
the six-vertex and the four-vertex theorems on plane curves.
Examples of quartic perturbations of the standard hyperboloid 
confirm our conjecture.
Our main result is a linearization and reformulation of the
problem in the framework of 2-dimensional Sturm theory;
we also define a signature of a quadratic point and calculate 
local normal forms recovering and
generalizing Tresse-Wilczynski's theorem.
\end{abstract}


\section{Introduction} \label{intro}

Almost one hundred years ago, S. Muchopadhyaya discovered two
theorems on plane ovals (an oval is a smooth closed strictly
convex plane curve). The first one is known as the four-vertex
theorem: the curvature of a plane oval has at least 4 critical
points. These critical points are the points at which the
osculating circles are hyperosculating, that is, are third-order
tangent to the curve. The second theorem concerns osculating
conics and states that a smooth convex closed curve has at least
6 distinct points at which the osculating conics are
hyperosculating. Such points are called
sextactic. A smooth plane curve can be approximated by
a conic at every point up to order 4; a point  is sextactic if
the order of approximation at this point is higher.  

The four- and six-vertex theorems and their ramifications continue
to attract interest, in great part due to work of V. I. Arnold who
placed the subject into the framework of symplectic and contact
topology \cite{ArnS,ArnU}. There is a wealth of new results in
this field, see \cite{O-T} for a survey.   

It is natural to expect that there exist multi-dimensional
versions of four- and six-vertex theorems but, so far, only the
very first steps have been made in this direction
\cite{Arn1,Arn2, Pa, Ur}. 

We consider the classical Carath\'eodory
conjecture as belonging to the area. This conjecture states that
a sufficiently smooth convex closed surface in $\R^3$ has at
least 2 distinct umbilic points, that is, the points where the two
principal curvatures are equal (see, e.g., \cite{O-T} and
references therein for a long and convoluted history of the
subject). Umbilic points are analogs of vertices of plane curves: these are the points
at which a sphere is abnormally (second-order) tangent to the surface. 
Let us note that a generic closed surface, even an immersed one, 
carries at least 4 umbilic points.

A smooth hypersurface $M$ in $\RP^3$ can be approximated by a quadric at
every point up to order 2. A point 
$x\in M$ is called \textit{quadratic} if $M$ can be approximated
by a quadric at $x$ up to order 3.\footnote{Quadratic points are
also called {\it hyperbonodes}, see \cite{Ur}; in \cite{Pa} these points are
called {\it special}.} We view quadratic points
of surfaces as 2-dimensional analogs of sextactic points. 
Quadratic points were studied in the
classical literature, see \cite{Wil, Sa, CF}, but
we are aware of only one existence result: 
if a generic smooth surface in $\RP^3$
contains a hyperbolic disc, bounded by a Jordan parabolic curve, then there exists
an odd number of quadratic points inside this disc (and hence, at least one)
\cite{Ur}.

\begin{figure}[ht]
\centerline{\epsfbox{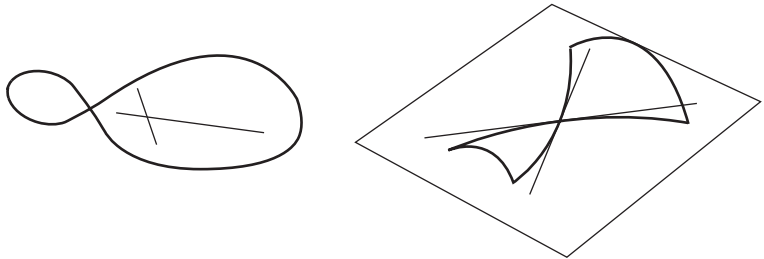}}
\caption{Non-degenerate hyperbolic surface and its asymptotic
directions at a generic point.}
\label{NonDe}
\end{figure}

We assume that $M$ is an orientable 
\textit{non-degenerate hyperbolic} surface: the second quadratic
form is non-degenerate and indefinite everywhere. Clearly, $M$ is
diffeomorphic to the 2-torus: $M\cong\T^2$. Indeed, at each point
of $M$ one has two \textit{asymptotic directions}, the light-cone
of the second quadratic form (equivalently defined as the tangent
lines to the intersection of $M$ with its tangent plane, see
figure \ref{NonDe}). It follows that the Euler characteristic of
$M$ is zero. The integral curves corresponding to the asymptotic
directions are called \textit{asymptotic lines} and they form a
2-web on $M$.
\begin{itemize}
\item {\it How many quadratic points are there on $M$?}
\end{itemize}

We start with local analysis and define a 
\textit{signature} of a non-gedenerate quadratic
point: $(s_1,s_2)$, with $s_i=\pm1$, this is a $\PGL(4,\R)$-invariant.
Explicit formul{\ae} for \textit{normal forms} provide further
differential invariants.
This problem goes back to Tresse and Wilczynski, we give here
a short proof of their classical result. We then calculate
normal forms at the degeneration stratas that have been studied
in \cite{Pl,land}.

An example of a non-degenerate hyperbolic surface is a
hyperboloid, $\cH$, given, in homogeneous coordinates, by the
equation 
\begin{equation}
\label{Hyper} x_0x_3=x_1x_2. 
\end{equation} 
Every hyperbolic quadric in $\RP^3$ is equivalent
to $\cH$ with respect to the action of the projective group
$\PGL(4,\R)$. 
Given a generic perturbation of $\cH$ given by a smooth
periodic function $f(u,v)$, we will prove that the quadratic points
are the points $(u,v)$ for which 
\begin{equation}
\label{SSyst}
\left\{
\begin{array}{rcl} 
f_{uuu}+f_{u}&=&0,\\[4pt] 
f_{vvv}+f_{v}&=&0.
\end{array}
\right.
\end{equation}

The Sturm-Hurwitz theorem states that a smooth
periodic function has no fewer zeroes than its first non-trivial
harmonic. In particular, the equation $f'''(x)+f'(x)=0$ has at
least four distinct roots on the circle $[0,2\pi)$ for every
$2\pi$-periodic function $f(x)$. This result implies the
classical four-vertex theorem.
In a recent paper \cite{P-Ch}  the following 
conjecture of V. Arnold is proved: if a plane wave front is Legendrian 
isotopic to a circle then it has at least four vertices. 
The vertices correspond 
to the solutions of the system
$$
\left\{
\begin{array}{rcl} 
F_{uuu}+F_{u}&=&0,\\[4pt] 
F_{v}&=&0
\end{array}
\right.
$$
where $F(u,v)$ is a generating function of the corresponding Legendrian curve 
in the space of cooriented 
contact elements of the plane (contactomorphic to the jet space $J^1 S^1$);
$u\in S^1$ is a cyclic coordinate and $v\in\R^k$ is an auxiliary variable.
One cannot help noticing that the above system bears a strong resemblance of 
our system (\ref{SSyst}); we believe that both are particular cases of 
multidimensional Sturm theory yet to be discovered. 

In general, we do not know how to estimate below
the number of solutions of (\ref{SSyst}). We will restrict
ourselves to the case where $f$ is a trigonometric polynomial of
bidegree $\leq(2,2)$ and prove a number of partial results. 
The geometric meaning of trigonometric polynomials of bidegree
$\leq(2,2)$ is that this class of functions $f$ describes the
perturbations of the hyperboloid that lie on a 
\textit{quartic}.
This way our considerations are related to an interesting problem of real
algebraic geometry: study quadratic points on quartics.\footnote{
For (non-degenerate)
cubic surfaces the situation is well understood. The
quadratic points in this case are precisely the intersection
points of the lines that lie on the surface (there are 27 complex
lines, but not all of them must be real), see \cite{Se,Pa}. In
particular, if a cubic surface is diffeomorphic to $\RP^2$, then
it has exactly 3 quadratic points \cite{Pa}.}
It is worth mentioning that the topology of a quartic that has a
$\T^2$-component $C^\infty$-close to the standard hyperboloid is
known. Such a quartic can have two components diffeomorphic to $\T^2$,
or one $\T^2$-component with $n$ spheres $S^2$, where $n=0,1,\ldots,9$,
see \cite{DK} for details.

Based on our partial results we conjecture that (\ref{SSyst}) has
at least 8 distinct zeros. This would imply that a small
perturbation of the hyperboloid has no less than 8 distinct
quadratic points. Let us make a bolder conjecture: \textit{every
closed hyperbolic surface in $\RP^3$ has no less than 8 distinct
quadratic points.} This conjecture is in the same relation to the
Carath\'eodory conjecture as the six-vertex theorem to the four-vertex one.

\section{Local analysis} \label{one}

In this section we formulate our problem and study
the local invariants of hyperbolic surfaces and quadratic points.

\subsection{Non-degenerate surface and quadratic points:\\
definitions and simple properties}

We collect here simple facts about quadratic points. Most
of them are known and can be found in classical books, see
\cite{Sa,CF,lane}.

Identify locally $\RP^3$ with the Euclidean space $\R^3$ with
coordinates 
\begin{equation}
\label{affCoord} x=x_1/x_0,
\quad y=x_2/x_0,
\quad z=x_3/x_0.
\end{equation} Given a hyperbolic surface $M$, these coordinates
can be chosen in such a way that in a neighbourhood of a point
$m$ this surface is given by
\begin{equation}
\label{LocCoord}
\textstyle z= xy+
\frac{1}{3}\left( a\,x^3+b\,y^3\right)+
\frac{1}{2}\left( c\,x^2y+d\,xy^2\right) +O(4)
\end{equation} 
where $a,b,c,d$ are some constants. Indeed, it
suffices to chose the asymptotic directions at $m$ as the 
$x$- and $y$-axes.

It is important to notice that the parameters $a,b,c,d$ are \textit{not}
well-defined functions of $m$. These parameters depend on the choice of
coordinates, for instance, the coordinate changes
$(x,y)\mapsto (tx,t^{-1}y)$ with arbitrary $t$ preserve the form of quation
(\ref{LocCoord}) but vary the parameters.
Even the signs of $a$ and $b$ change as one changes $(x,y)$ to $(-x,-y)$.
The geometric meaning of $a$ and $b$ will be explained in Appendix, see also
\cite{O-T}.

Nevertheless, zero sets $a=0$ and $b=0$ are well defined. 

\begin{fact}
\label{fact-1} 
Point $m$ is a quadratic point if and only if the paramaters $a$ and $b$ in
(\ref{LocCoord}) vanish at $m$:
\begin{equation}
\label{CondQuad}
\left\{
\begin{array}{rcl} a&=&0\\
b&=&0.
\end{array}
\right.
\end{equation}
\end{fact}

\begin{proof}
First we check that the condition (\ref{CondQuad}) is independent of the choice of
coordinates $x$ and $y$.

An osculating quadrics at a point $m$ is as follows
\begin{equation}
\label{OscQuad}
\textstyle z=xy+\frac{1}{2}\left(
\gamma\,xz+\delta\,yz+\eps\,z^2
\right),
\end{equation} 
where $\gamma,\delta$ and $\eps$ are arbitrary constants.
Indeed, formula (\ref{OscQuad}) defines the quadrics
approximating $M$ up to the terms of order 2.
Let now $m$ be a quadratic point. A quadric (\ref{OscQuad}) is
\textit{hyperosculating} if it coincides with $M$ up to the terms
of order 3. This is the case if and only if $\gamma=c,\delta=d$ and
the constants $a$ and $b$ in (\ref{LocCoord}) vanish. \proofend
\end{proof}

Let us summarize the above calculations.
\begin{fact}
\label{fact0} 
(i)
At a generic point, there exists a 3-parameter family of osculating
quadrics given by formula (\ref{OscQuad}).

(ii)
At a quadratic point, there is a 1-parameter family of
hyperosculating quadrics. 
\end{fact}

\noindent
Indeed, $\eps$ in (\ref{OscQuad}) remains a free parameter.

We can now define explicitly the notion of a generic surface that will be
essential for the sequel. 
The following definition is open and dense in $C^\infty$-topology.

\begin{definition}
{\rm
A non-degenerate hyperbolic surface $M$ in $\RP^3$ is said to be generic, or in
general position, if:

\noindent
1) the sets $(a=0)$ and $(b=0)$ are smooth embedded curves in $M$ with transversal
intersections;

\noindent
2) at each intersection point $(a=0)\cap(b=0)$ both curves $(a=0)$ and $(b=0)$ are
transversal to the asimptotic directions.
}
\end{definition}

We arrive at the following observation that justifies the
formulation of our main problem.

\begin{fact}
\label{fact1} Quadratic points on a generic hyperbolic surface in
$\RP^3$ are isolated.
\end{fact}
\noindent 
Indeed, condition (\ref{CondQuad}) is of
codimension 2 since the parameters $a$ and $b$ are two
independent functions in $(x,y)$.

A hyperbolic surface is
a quadric if and only if it contains its asymptotic tangent lines
at any point, cf. \cite{WilBis}. The next statement is nothing else but an
infinitesimal version of this statement (see, e.g., \cite{CF},
p. 62). 

\begin{fact}
\label{fact2} Quadratic points are those points at which both
asymptotic lines have inflections.
\end{fact}

The curves $(a=0)$ and $(b=0)$ on $M$ are precisely the sets of
inflection points of the two asymptotic foliations; the quadratic
points are the intersection points $(a=0)\cap(b=0)$.\footnote{The
union $(a=0)\cup(b=0)$ is usually called the {\it flecnodal
curve}.}

\begin{remark} 
{\rm For an arbitrary oriented smooth foliation on
$\T^2$, the average curvature of the leaves with respect to the
standard flat metric is zero, see \cite{As}. Therefore the leaves of any
foliation have inflection points. It is easy to find two transversal foliations
with no points at which both leaves have inflections, see figure
\ref{CounterExFig}.
\begin{figure}[ht]
\centerline{\epsfbox{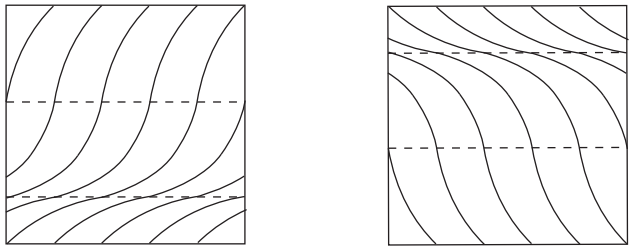}}
\caption{Transversal foliations with no common inflections.}
\label{CounterExFig}
\end{figure}
This would be a counterexample to our conjecture if
one could realise these foliations as asymptotic lines on a
hyperbolic surface. 
}
\end{remark}

\subsection{Signature of a quadratic point}

Recall that we consider only non-degenerate quadratic points.
In this section we define an invariant of such a quadratic point that we call \textit{signature}.

\begin{figure}[ht]
\centerline{\epsfbox{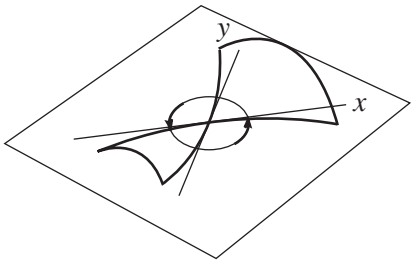}}
\caption{Natural ordering of $x$- and $y$-axis.}
\label{OrdFig}
\end{figure}

Fix an orientation of $M$ and of $\RP^3$; the surface $M$ is then
co-oriented.
The asymptotic $x$- and $y$- directions are
naturally \textit{ordered} at any point $m$.
Indeed, choose $z$-coordinate in (\ref{LocCoord}) positively
coorienting $M$, consider the tangent plane $T_mM$ and draw a small circle on it
centered at $m$. Choose a point on the circle which lies above $M$ and move in
the positive direction; the first intersection with $M$ corresponds to the
$x$-axis, see figure \ref{OrdFig}. 

\begin{figure}[ht]
\centerline{\epsfbox{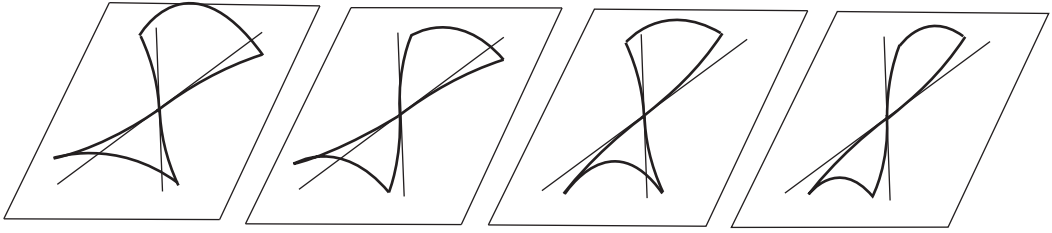}}
\caption{Quadratic points of signature 
$(+,+)$, $(+,-)$, $(-,+)$ and $(-,-)$.}
\label{SignFigLab}
\end{figure}

 At a quadratic point $m$, the surface $M$ can
``cross'' the tangent plane in four different ways, see figure
\ref{SignFigLab}. 

\begin{definition}
\label{signDef}
{\rm 
Define the signature $s=(s_1,s_2)$, where
$s_i=+$ or $-$, of a quadratic point $m$. We put $s_1=+$ (resp.
$s_2=+$) if the $x$-axis (resp. $y$-axis) in the vicinity of $m$
lies under $M$. We put the $-$ sign otherwise. 
}
\end{definition}

Clearly, the signature is a $\PGL(4,\R)$-invariant of a quadratic
point. Note that if one changes the orientation of $M$ or
$\RP^3$, then the signs $s_1$ and $s_2$ change:
$(+,+)\leftrightarrow(-,-)$ and $(+,-)\leftrightarrow(-,+)$. One
can call the points of the two above types \textit{even} and
\textit{odd}, respectively. This notion of parity is independent of
the choice of orientation.

For every quadratic point $m$, coefficients $a$ and $b$ from
(\ref{LocCoord}) vanish at $m$. Consider the expansion
(\ref{LocCoord}) for a point close to $m$.

\begin{lemma}
\label{ProstLem} 
One has $s_1=+$ if and only if $ax\geq 0$ on the $x$-axis; and
$s_2=+$ if and only if $by\geq0$ on the  $y$-axis.
\end{lemma}

\begin{proof} 
First notice that the coordinate change $(x,y)\mapsto(-x,-y)$ changes the signs
of $a$ and $x$ simultaneously (as well as signs of $b$ and $y$),
so that the signs of the expressions $a\,x$ and $b\,y$ are well defined. 

Consider the expansion (\ref{LocCoord}) on the positive $x$-semiaxis.
Since $y=0$, one has
$z(x,0)=\frac{1}{3}\,a\,x^3+O(4)$. By definition of signature, $s_1=+$ means that
$z(x,0)>0$. The curve $(a=0)$ is
transversal to the $x$-axis, and the statement follows. \proofend
\end{proof}

A family of non-degenerate hyperbolic surfaces $M_t$ smoothly
depending on a parameter $t\in[0,1]$ is called a
\textit{homotopy}. 
We do not assume \textit{a-priori} that at each moment $t$ the surface $M_t$
is generic.

\begin{figure}[ht]
\centerline{\epsfbox{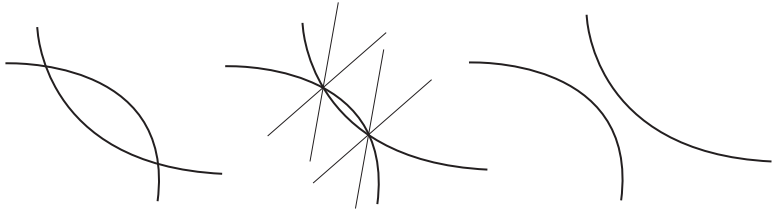}}
\caption{Creation/annihilation of quadratic ponits.}
\label{CreAnnFig}
\end{figure}

Quadratic points can be ``created'' or
``annihilated'' by homotopy in pairs, see
figure~\ref{CreAnnFig}.

\begin{proposition}
\label{CreAnnProp} 
Two quadratic points created/annihilated by a
homotopy are of the same signature.
\end{proposition}
\begin{proof} Close to the moment of creation/annihilation of a
pair of quadratic points the curves $(a=0)$ and $(b=0)$ are
transversal to both asymptotic directions. The statement then
follows from Lemma \ref{ProstLem}. \proofend
\end{proof}

Consider now a homotopy in the class of generic surfaces, i.e.,
$M_t$ is generic for all $t\in[0,1]$. Let us call such a homotopy
\textit{stable}. Signature is preserved by a
stable homotopy. 

\begin{figure}[ht]
\centerline{\epsfbox{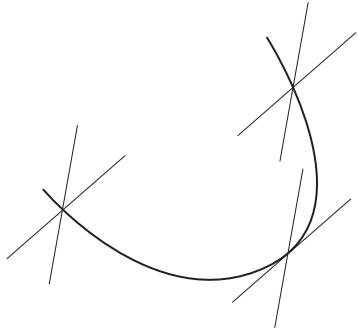}}
\caption{Curves $(a=0)$ becomes non-transversal to the $x$-axis.}
\label{NonTransFig}
\end{figure}

\begin{lemma}
\label{StabProp} If a quadratic point $m\in{}M$ is not
annihilated by a stable homotopy $M_t$, then the signature of
$m_t$ does not depent on $t$.
\end{lemma}
\begin{proof} Suppose that the homotopy connects points $m_1$ and
$m_2$ of signature $(+,+)$ and $(-,+)$, respectively. Then there
is a moment $t_0$ at which the curve $(a=0)$ is not transversal
to the $x$-axis, see figure \ref{NonTransFig}. 
This contradicts to the fact that $M_{t_0}$ is in general position. \proofend
\end{proof}

\subsection{Normal forms and differential invariants} \label{NFSec}

The normal form of a non-degenerate
hyperbolic surface $M$ in the vicinity of a generic point $m$ is one
of the most classical results of projective differential geometry
that goes back to Tresse and Wilczynski.
In this section we give a much simpler proof of Tresse-Wilczynski's
result; we then calculate the normal form (up to the order 4) 
in the vicinity of a point that lies
on the curve $(a=0)$ but not on $(b=0)$, 
and finally in the vicinity of a quadratic point.

Normal forms are also discussed in \cite{Pl} and \cite{land}, but
the formula $\Pi_{3,1}$ of the former paper differs from the Tresse-Wilczynski's result,
while the formula $\Pi_{4,3}$ is different from our formula
(\ref{NewNFEq}) below.

\begin{theorem}
\label{norforThm}
Modulo projective transformations, a non-degenerate hyperbolic surface $M$
is given, in a vicinity of a point $m$, by the formulas:

(i)
if $m$ is generic then
\begin{equation}
\label{WilNFEq}
\textstyle 
z= xy+
\frac{1}{3}\left(
x^3+y^3
\right)+
\frac{1}{12}\left(
I\,x^4+J\,y^4
\right)
+O(5),
\end{equation}

(ii)
if $m$ belongs to the curve $(a=0)$ but not to $(b=0)$ then
\begin{equation}
\label{IntNFEq}
\textstyle 
z= xy+
\frac{1}{3}\left(
y^3\pm
x^3y
\right)+
\frac{1}{12}\,
\tilde{I}\,x^4
+O(5),
\end{equation}

(iii)
if $m$ is quadratic then
\begin{equation}
\label{NewNFEq}
\textstyle 
z= xy\pm
\frac{1}{3}\left(
x^3\,y\pm x\,y^3
\right) 
+\frac{1}{12}\left(
\bar{I}\,x^4+\bar{J}\,y^4
\right)
+O(5),
\end{equation}
(all four combinations of signs are possible),
where the parameters $(I,J)$ are $\PGL(4,\R)$-invariants, as well as
$\tilde{I}$ defined up to the sign; $(\bar{I},\bar{J})$ are invariants defined
up to the simultaneous sing change.
If $M$ is oriented then the signs of $\tilde{I}$ and $(\bar{I},\bar{J})$
are well-defined.
\end{theorem}

\begin{proof}
Consider the action of the Lie algebra $\Sl(4,\R)$ which is the infinitesimal
version of the $\PGL(4,\R)$-action. In affine coordinates $(x,y,z)$, this
action is spanned by 3 constant vector fields of ``translations'', together with 9
linear vector fields and 3 quadratic vector fields of ``inversions'':
\begin{equation}
\label{Sl4act}
\textstyle 
\Sl(4,\R)=
\left\langle
\frac{\partial}{\partial x},
\frac{\partial}{\partial y},
\frac{\partial}{\partial z},
\quad
x\,\frac{\partial}{\partial x},
\ldots,
z\,\frac{\partial}{\partial z},
\quad
x\,\cE,
y\,\cE,
z\,\cE
\right\rangle,
\end{equation}
where 
$
\textstyle 
\cE=
x\,\frac{\partial}{\partial x}+
y\,\frac{\partial}{\partial y}+
z\,\frac{\partial}{\partial z}
$
is the Euler field.

Consider a hyperbolic surface $M$
$$
z=xy+O(3)
$$
where $O(3)$ stands for functions of $x$ and $y$ that belong to the cube of the maximal ideal $(x,y)$. 
One readily checks the following
\begin{lemma}
The Lie algebra of vector fields preserving the second jet of $M$ is the 
subalgebra of dimension 7 spanned by
\begin{equation}
\label{subalgEq}
\textstyle 
\left\langle
x\frac{\partial}{\partial x}
+z\frac{\partial}{\partial z},\,
y\frac{\partial}{\partial y}+
z\frac{\partial}{\partial z},
\quad
z\frac{\partial}{\partial x},\,
z\frac{\partial}{\partial y},\,
\quad
x\,\cE,\,
y\,\cE,\,
z\,\cE
\right\rangle.
\end{equation}
\end{lemma}
Consider the expansion (\ref{LocCoord}) and let us study the action of the Lie algebra
(\ref{subalgEq}) on the coefficients $(a,b,c,d)$. 
The action of $X=\l\,x\,\cE+\m\,y\,\cE$ is:
$$
\dot{a}=0,
\quad
\dot{b}=0,
\quad
\dot{c}=2\l,
\quad
\dot{d}=2\m,
$$
where $\dot{}$ stands for the Lie derivative $L_X$; it follows that the flow of
such an element ``kills'' the coefficients $c$ and $d$. 

The action of
$X=\nu\,(x\frac{\partial}{\partial x}
+z\frac{\partial}{\partial z})+\kappa\,(y\frac{\partial}{\partial y}+
z\frac{\partial}{\partial z})$ is:
$$
\dot{a}=3\nu a,
\quad
\dot{b}=3\kappa b,
\quad
\dot{c}=2\nu c,
\quad
\dot{d}=2\kappa d.
$$
One concludes that the expansion of $M$ can be reduced to the form
\begin{equation}
\label{OrTreGen}
\textstyle
z= xy+
\frac{1}{3}\left(
x^3+y^3
\right)+O(4)
\end{equation}
if $m$ is generic so that $a\not=0$, $b\not=0$; to
\begin{equation}
\label{OrTreCDOne}
\textstyle
z= xy+
\frac{1}{3}\,
y^3+O(4),
\end{equation}
if $m$ belongs to the curve $(a=0)$ but with $b\not=0$, and to
\begin{equation}
\label{OrTreQua}
z= xy+O(4),
\end{equation}
if $m$ is quadratic.
Indeed, one can assume, without loss of generality, that
$a\geq0$ and $b\geq0$ (it suffices to change the coordinates $(x,y,z)$ to $(-x,y,-z)$, 
or $(x,-y,-z)$, or $(-x,-y,z)$ to change the signs of $a$ and $b$).
One then finds a vector field from
(\ref{subalgEq}) whose flow reduces the coefficients to $a,b$ to $1$; whenever
they are different from $0$.

Part (i).
If $m$ is generic then
the subalgebra of (\ref{subalgEq}) preserving the
third-order expansion (\ref{OrTreGen}) is of dimensionl 3 and spanned by
$$
\textstyle
\langle
z\,\frac{\partial}{\partial x}
-y\,\cE,\,
z\,\frac{\partial}{\partial y}
-x\,\cE,\,
z\,\cE
\rangle,
$$
Consider an arbitrary 4-th order expression
\begin{equation}
\label{ArbQua}
Q_4(x,y)=\a\,x^4+\b\,x^3y+\g\,x^2y^2+\de\,xy^3+\eps\,y^4,
\end{equation}
the action of 
$X=\l\,(z\,\frac{\partial}{\partial x}-y\,\cE)+
\m\,(z\,\frac{\partial}{\partial y}-x\,\cE)+
\nu\,z\,\cE$ 
is given by
$$
\textstyle
\dot{\a}=-\frac{1}{3}\,\m,
\quad
\dot{\b}=\frac{2}{3}\,\m,
\quad
\dot{\g}=\nu,
\quad
\dot{\de}=\frac{2}{3}\,\l,
\quad
\dot{\eps}=-\frac{1}{3}\,\l,
$$
so that one can kill the coefficients $\b,\g$ and $\de$.
It follows that (\ref{WilNFEq}) is the normal form of $M$ in a neighbourhood of a generic point.

Part (ii).
The subalgebra of (\ref{subalgEq}) preserving (\ref{OrTreCDOne}) is spanned by 4 vector
fields
$$
\textstyle
\langle
\frac{2}{3}\,x\frac{\partial}{\partial x}
+\frac{1}{3}\,y\frac{\partial}{\partial y}
+z\frac{\partial}{\partial z},\;
z\,\frac{\partial}{\partial x}
-y\,\cE,\;
z\,\frac{\partial}{\partial y}
-x\,\cE,\;
z\,\cE
\rangle.
$$
As above, the action of $X=z\,\cE$ allows one to kill the coefficient $\g$ in (\ref{ArbQua}).
The action of 
$X=\l(z\,\frac{\partial}{\partial x}-y\,\cE)+
\m(z\,\frac{\partial}{\partial y}-x\,\cE)$
on the 4-th order part reads
$$
\textstyle
\dot{\a}=0,
\quad
\dot{\b}=0,
\quad
\dot{\g}=0,
\quad
\dot{\de}=\frac{2}{3}\,\m,
\quad
\dot{\eps}=-\frac{1}{3}\,\l
$$
that kills the coefficients $\de$ and $\eps$ in (\ref{ArbQua}).
Finally, the action of the vector field
$X=\frac{2}{3}\,x\frac{\partial}{\partial x}
+\frac{1}{3}\,y\frac{\partial}{\partial y}
+z\frac{\partial}{\partial z}$ is
$$
\textstyle
\dot{\a}=\frac{5}{3}\,\a,
\quad
\dot{\b}=\frac{4}{3}\,\b,
\quad
\dot{\g}=\g,
\quad
\dot{\de}=\frac{2}{3}\,\de,
\quad
\dot{\eps}=\frac{1}{3}\,\eps,
$$
so that the coefficient $\b$ can be reduced to $\pm1$.
Formula (\ref{IntNFEq}) is proved.
Simultaneous change of the signs: $(y,z)\leftrightarrow(-y,-z)$
changes the sign of $\tilde{I}$. 
This, of course, changes the (co)orientation of $M$ defined by
the $z$-axis.

Part (iii).
The subalgebra of (\ref{subalgEq}) preserving the
third-order expansion (\ref{OrTreQua}) is spanned by 5 vector
fields
$$
\textstyle
\langle
x\frac{\partial}{\partial x}
+z\frac{\partial}{\partial z},\;
y\frac{\partial}{\partial y}+
z\frac{\partial}{\partial z},\;
z\,\frac{\partial}{\partial x}
-y\,\cE,\;
z\,\frac{\partial}{\partial y}
-x\,\cE,\;
z\,\cE
\rangle.
$$
The action of 
$X=\l\,(x\,\frac{\partial}{\partial x}+z\,\frac{\partial}{\partial z})+
\m\,(y\,\frac{\partial}{\partial y}+z\,\frac{\partial}{\partial z})$
is
$$
\textstyle
\dot{\a}=(3\l-\m)\,\a,
\quad
\dot{\b}=2\l\b,
\quad
\dot{\g}=(\l+\m)\,\g,
\quad
\dot{\de}=2\m\de,
\quad
\dot{\eps}=(3\m-\l)\,\eps,
$$
so that one can reduce $\b$ and $\de$ to $\pm1$ if only these coefficients are 
different from zero, and this is the case since the quadratic point is generic . 
As above, one reduces the coefficient $\g$ in (\ref{ArbQua}) to 0. 
The actions of the fields 
$z\,\frac{\partial}{\partial x}-y\,\cE$ and 
$z\,\frac{\partial}{\partial y}-x\,\cE$ is trivial. 
Formula (\ref{NewNFEq}) follows.

Again, changing the signs: $(x,z)\leftrightarrow(-x,-z)$ 
or $(y,z)\leftrightarrow(-y,-z)$, 
one changes the signs $(\bar{I},\bar{J})\leftrightarrow(-\bar{I},-\bar{J})$.
This simultaneous sign change
corresponds to the change of the orientation.

Theorem \ref{norforThm} is proved. \proofend
\end{proof}

\begin{remark}
{\rm
Formula (\ref{WilNFEq}) is precisely the normal form of Tresse and Wilczynski 
(see \cite{Wil}, Second Memoir, formula (96)).
The coefficients $I,J$ and all of the following coefficients are called 
\textit{absolute invariants}\footnote{In Fifth
Memoir Wilczynski develop the series up to order 6 and interpret the next
13 coefficients.} of $M$.
}
\end{remark}

\begin{lemma}
\label{SInvLem}
The signature of a quadratic point is nothing else but the sign of the invariants
in (\ref{NewNFEq}), namely
$(\sigma_1,\sigma_2)=(\mathrm{sgn}\bar{I},\mathrm{sgn}\bar{J})$.
\end{lemma}

\begin{proof}
This follows directly from Definition \ref{signDef}. 
Indeed, restricting the right hand side of (\ref{NewNFEq}) to the $x$-axis, one has:
$z=\bar{I}\,x^4+O(5)$. Hence $\bar{I}>0$ if and only if $M$ is above the $x$-axis, and likewise for the $y$-axis. \proofend
\end{proof}

\section{Small perturbations of the hyperboloid:\\ linearization
of the problem} \label{two}

In this section we deduce system (\ref{SSyst}) as the  first-order
approximation to our problem.

A \textit{perturbation} of the hyperboloid $\cH$ is a homotopy
$M_\eps$, smoothly depending on a small parameter
$\eps\in\R$, such that $M_0=\cH$. When we talk of ``sufficiently
small'' perturbations, this means that there exists $\eps_0>0$ such
that the property we consider holds for all
$|\eps|\leq\eps_0$. 

The hyperboloid $\cH$ defined by formula (\ref{Hyper}) has the
following natural parametrization:
\begin{equation}
\label{Param}
\begin{array}{rcl} x_0(u,v)&=&\cos\frac{u}{2}\,
\cos\frac{v}{2},\\[4pt] 
x_1(u,v)&=&
\cos\frac{u}{2}\,
\sin\frac{v}{2},\\[4pt] 
x_2(u,v)&=&
\sin\frac{u}{2}\,
\cos\frac{v}{2},\\[4pt] 
x_3(u,v)&=&
\sin\frac{u}{2}\,
\sin\frac{v}{2},
\end{array}
\end{equation} where $(u,v)\in[0,2\pi)$. The coordinates $(u,v)$
on $\cH$ are globally defined. Although $x_i(u,v)$ are
not well defined functions on the torus, formula (\ref{Param})
gives a well-defined embedding $\T^2\hookrightarrow\RP^3$. 

We can describe a small perturbation of $\cH$ in terms of a
function on $\T^2$; the construction is as follows. The ``normal''
vector
$X_{uv}:=\frac{\partial^2}{\partial{}u\partial{}v}\,X(u,v)$, given
more explicitly by
\begin{equation}
\label{ProjNor}
X_{uv}
=
\textstyle
\frac{1}{4}\,
\left(
\sin\frac{u}{2}\,\sin\frac{v}{2},\quad
-\sin\frac{u}{2}\,\cos\frac{v}{2},\quad
-\cos\frac{u}{2}\,\sin\frac{v}{2},\quad
\cos\frac{u}{2}\,\cos\frac{v}{2}
\right),
\end{equation} 
is always transversal to $\cH$, and so the family of surfaces
\begin{equation}
\label{PerFunF}
\widetilde{X}(u,v)=X(u,v)+\eps\,f(u,v)\,X_{uv},
\end{equation} 
where $f:\T^2\to\R$ is an arbitrary smooth function, remains
smooth for sufficiently small $\eps$. Conversely, every surface
$M$ sufficiently close to $\cH$ can be represented in a
parametrized form by (\ref{PerFunF}).

\begin{proposition}
\label{SElProp} 
The perturbed surface (\ref{PerFunF}) remains a
quadric, in the first order in $\eps$, if and only if the
function $f$ is a combination of the first harmonics:
\begin{equation}
\label{FirHar} 
f=\sum_{-1\leq{}n,m\leq{}1}
f_{n,m}\,e^{i(n\,u+m\,v)},
\end{equation} 
where $f_{n,m}\in \C$ and
$f_{n,m}=\overline{f}_{-n,-m}$ (since $f$ is real).
\end{proposition}

\begin{proof} From (\ref{Param}) and (\ref{ProjNor}) one readily obtains the
equation of the perturbed surface (\ref{PerFunF}):
\begin{equation}
\label{PerturFEq}
\widetilde{x}_0\widetilde{x}_3 
-\widetilde{x}_1\widetilde{x}_2
=\frac{\eps}{4}\,f+O(\eps^2).
\end{equation} 
We will need the following

\begin{lemma}
\label{HarQuadLem}
Function $f$ is a combination of the first harmonics if and only
if $f$ is a quadratic expression in the coordinates $x_i(u,v)$
given by (\ref{Param}):
$$ 
f(u,v)=\sum\alpha_{ij}\,x_i\,x_j,
$$ 
where $\alpha_{ij}$ are arbitrary constants.
\end{lemma}
\begin{proof}
The proof of this lemma is quite obvious. For instance, one has
$$
\textstyle
\sin{}u\,\sin{}v=4\,\cos\frac{u}{2}\,
\sin\frac{u}{2}\,\cos\frac{v}{2}\,
\sin\frac{v}{2}=4\,x_0\,x_4
$$
and similarly for other homogeneous first-order harmonics,
whereas
$$
\textstyle
\sin{}u=2\,\cos\frac{u}{2}\,
\sin\frac{u}{2}
\left(
\cos^2\frac{v}{2}+
\sin^2\frac{v}{2}
\right)=2\left(
x_0\,x_2+x_1\,x_4
\right)
$$
and similarly for other homogeneous harmonics of order $(0,1)$ or
$(1,0)$, and finally one has
$$
1=x_0^2+x_1^2+x_2^2+x_3^2+x_4^2
$$
for the constant function.
Lemma \ref{HarQuadLem} follows. \proofend
\end{proof}

Equation (\ref{PerturFEq}) implies now that the perturbed surface
(\ref{PerFunF}) satisfies the quadratic equation:
$$
\widetilde{x}_0\widetilde{x}_3 
-\widetilde{x}_1\widetilde{x}_2
-\frac{\eps}{4}\,
\sum_{0\leq{}i,j\leq{}3}
\alpha_{ij}\,\widetilde{x}_i\widetilde{x}_j +O(\eps^2)=0.
$$ 
Therefore, the perturbed surface (\ref{PerFunF}) remains
quadric in the first order in $\eps$.
Proposition \ref{SElProp} is proved. \proofend
\end{proof}

Notice that the space of functions (\ref{FirHar}) is precisely
the 9-dimensional space of solutions of the system
$$
\left\{
\begin{array}{rcl} f_{uuu}+f_{u}&\equiv &0,\\[6pt]
f_{vvv}+f_{v}&\equiv &0.
\end{array}
\right.
$$ The following statement can be understood as a version of
Proposition \ref{SElProp}, but at a single point.

\begin{theorem}
\label{SProp} Given a perturbation (\ref{PerFunF}) of the
standard hyperboloid, a point with coordinates $(u,v)=(u_0,v_0)$
remains quadratic in the first order in $\eps$ if and only if
condition (\ref{SSyst}) is satisfied at $(u_0,v_0)$.
\end{theorem}

\begin{proof} 
Without loss of gererality consider the point with
coordinates $(u_0,v_0)=(0,0)$. Identify locally $\RP^3$ and the
Euclidean space $\R^3$ with the coordinates (\ref{affCoord}); the
parametrized hyperboloid $\cH$ is then given by
\begin{equation}
\label{CoTan}
\begin{array}{rcl} x&=&\tan\frac{v}{2},\\[4pt]
y&=&\tan\frac{u}{2},\\[4pt] z&=&\tan\frac{u}{2}\,\tan\frac{v}{2}.
\end{array}
\end{equation} 
Let us calculate the perturbed surface
(\ref{PerFunF}). One obtains
$$
\widetilde{x}=
\frac{\cos\frac{u}{2}\,\sin\frac{v}{2}-\frac{\eps}{4}\,f\,\sin\frac{u}{2}\,\cos\frac{v}{2}}
{\cos\frac{u}{2}\,\cos\frac{v}{2}+
\frac{\eps}{4}\,f\,\sin\frac{u}{2}\,\sin\frac{v}{2}}.
$$ and similarly for $\widetilde{y}$ and $\widetilde{z}$. Finally
one has
\begin{equation}
\label{TilDeq}
\begin{array}{rcl}
\widetilde{x}&=& x-\frac{\eps}{4}\,f\left( y+xz\right) +O(\eps^2)
\\[6pt]
\widetilde{y}&=&y-\frac{\eps}{4}\,f \left( x+yz\right) +O(\eps^2)
\\[6pt]
\widetilde{z}&=&xy+\frac{\eps}{4}\,f \left( 1- z^2\right)
+O(\eps^2).
\end{array}
\end{equation}

Assume that the point
$(\widetilde{x},\widetilde{y},\widetilde{z})=(0,0,0)$ of the
perturbed surface is quadratic. This means its coordinates have
to satisfy a quadratic equation up to the terms of order $\leq2$
in $\eps$ and $\leq4$ in
$(\widetilde{x},\widetilde{y},\widetilde{z})$, namely
$$
\widetilde{z} -
\widetilde{x}\widetilde{y}=
\eps\left(
P(\widetilde{x},\widetilde{y},\widetilde{z})+O(4)\right)+
O(\eps^2),
$$ where $P$ is a polynomial of degree $\leq2$.

Exactly as in the case of condition (\ref{CondQuad}), the
coefficients of $x^3$ and $y^3$ in the above equation are the
obstructions to existence of such a polynomial $P$. Indeed, these
coefficients are identically zero (up to order 1 in $\eps$) in
the right hand side. Let us calculate these coefficients for the
left hand side of the above equality.

From the Taylor expansion one obtains, for the case of $x^3$, the following
expression:
$\left(\frac{1}{24}\,f_{xxx}+\frac{1}{4}\,f_x\right)\eps$, where
the derivatives are takes at the point $(0,0)$. In the same way,
one gets
$\left(\frac{1}{24}\,f_{yyy}+\frac{1}{4}\,f_y\right)\eps$ for
$y^3$. Therefore, one obtains the following system:
$$
\left\{
\begin{array}{rcl}
\frac{1}{6}\,f_{xxx}+f_x &=&0\\[6pt]
\frac{1}{6}\,f_{yyy}+f_y &=&0.
\end{array}
\right.
$$ Furthermore, the chain rule, applied to (\ref{CoTan}), implies that at point $(0,0)$ one has:
$$ f_x=2\,f_v
\quad
\hbox{and}
\quad f_{xxx}=8\,f_{vvv}-4\,f_{v},
$$ so that the above system is precisely the system (\ref{SSyst}). \proofend
\end{proof}

The following statement is an immediate consequence of Theorem
\ref{SProp} and of compactness of the 2-torus.

\begin{corollary} 
Given a deformation $M_\eps$ defined by
(\ref{PerFunF}) with sufficiently small $\eps$ and a generic function $f$, the number of
quadratic points on $M_\eps$ coincides with the number of the
points for which the system (\ref{SSyst}) is satisfied.
\end{corollary}

Indeed, for a generic function $f$, the solution of (\ref{SSyst}) are simple
(of multiplicity 1) and cannot be removed by a small perturbation.

\section{Approximation by quartics:
second harmonics} \label{three}

We will be interested in the perturbations of the hyperboloid
$\cH$ in the class of quartics. More precisely, we will be
looking for $C^\infty$-families $M_\eps$ of quartics that contain
a smooth component diffeomorphic to $\T^2$ and coinciding with
$\cH$ for $\eps=0$.

According to Proposition \ref{SElProp}, the
space of first harmonics (\ref{FirHar}) corresponds to the
perturbations of $\cH$ inside the space of
quadrics. It turns out that the space of second harmonics also
has a nice algebraic geometry meaning. 

\begin{proposition}
\label{QuarDefThm} 
A perturbation (\ref{PerFunF}) satisfies a
quartic equation, in the first order in $\eps$, if and only if
the function $f$ is given by the formula
\begin{equation}
\label{SecHar} 
f=\sum_{-2\leq{}n,m\leq{}2} f_{n,m}\,e^{i
(n\,u+m\,v)},
\end{equation} 
where $f_{n,m}\in \C$ and
$f_{n,m}=\overline{f}_{-n,-m}$.
\end{proposition}

\begin{proof} 
Function $f$ is as in (\ref{SecHar}) if and only if
$f$ can be written in terms of the coordinates (\ref{Param}) as a
homogeneous quartic expression:
$$ f=\sum_{0\leq{}i,j,k,\ell\leq{}3}
\alpha_{ijk\ell}\, x_ix_jx_kx_\ell
$$ where $\alpha_{ijk\ell}$ are some constants.
The proof of this statement is similar to that of Lemma
\ref{HarQuadLem}.

Equation (\ref{PerturFEq}) implies then that the perturbed surface
(\ref{PerFunF}) satisfies a homogeneous equation of order~4
$$
\left(\widetilde{x}_0\widetilde{x}_3
-\widetilde{x}_1\widetilde{x}_2\right)
\left(\widetilde{x}_0^2+\widetilde{x}_1^2+
\widetilde{x}_2^2+\widetilde{x}_3^2\right) 
-\frac{\eps}{4}\,
\sum_{0\leq{}i,j,k,\ell\leq{}3}
\alpha_{ijk\ell}\,
\widetilde{x}_i\widetilde{x}_j\widetilde{x}_k\widetilde{x}_\ell
+O(\eps^2)=0,
$$ 
and Proposition \ref{QuarDefThm} follows. \proofend
\end{proof}

Let us calculate the dimension of the moduli space of quartic
deformations of $\cH$.

\begin{proposition} 
The space of $\PGL(4,\R)$-classes of quartic
deformations of $\cH$ is of dimension 15.
\end{proposition}

\begin{proof} 
We give two ways to calculate the dimension of
the space of deformations.

\textit{First}. 
The space of second harmonics (\ref{SecHar}) is
25-dimensional. Its quotient by the space of first harmonics
(that do not change $\cH$ up to projective transformations, cf.
Proposition \ref{SElProp}) is 16. Finally, the quotient by
homotheties $\R^*$ leaves us with a 15-dimensional space.

\textit{Second}. 
The space $\R_4[x_0,x_1,x_2,x_3]$ of homogeneous
polynomials of degree 4 is of dimension 35 (and so the dimension
of the space of quartics is 34). The space of quartic
deformations modulo the
$\PGL(4,\R)$-action is related to the quotient space
$\R_4[x_0,x_1,x_2,x_3]/\mathcal{R}$, where
$\mathcal{R}$ is the component of degree 4 of the ideal with two generators:
$$
\mathcal{R}=
\left\langle 
x_0x_3-x_1x_2,
\quad x_0^2+x_1^2+x_2^2+x_3^2
\right\rangle.
$$ 
One easily checks that $\dim\mathcal{R}=19$, so that, taking
into account the homotheties, we again obtain dimension 15. \proofend
\end{proof}

Naturally, both calculation yield the same answer, in
accordance with Proposition \ref{QuarDefThm}. 

\section{Partial solutions to the main system} \label{four}

Unfortunately, we are unable to give a general estimate below of
the number of solutions of system (\ref{SSyst}). We will consider
the case where the function $f$ belongs to the space of second
harmonics, but even in this case our estimates are not complete.
We will give two partial results and one example that we believe
realizes the least number of solutions.

\medskip

{\bf A}.
Consider first a 12-dimensional subspace of the space
(\ref{SecHar}) with the condition
$f_{2,2}=f_{2,-2}=f_{-2,2}=f_{-2,-2}=0$, that is, the subspace 
of functions which are at most first harmonics in one of the
variables. 
\begin{proposition}
\label{OneTwo} If $f$ is a generic function belonging to the
above subspace, then there are at least 12 distinct points on the
$[0,2\pi)\times[0,2\pi)$-torus at which system (\ref{SSyst}) is
satisfied.
\end{proposition}
\begin{proof} One has:
$$
\begin{array}{rcl} f_{uuu}+f_u&=&
\phi_1(v)\,\cos2u+\phi_2(v)\,\sin2u,\\[6pt] f_{vvv}+f_v&=&
\psi_1(u)\,\cos2v+\psi_2(u)\,\sin2v
\end{array}
$$ where the functions $\phi_i$ and $\psi_i$ belong to the space
of first harmonics.

Consider first the curve $f_{vvv}+f_v=0$ on $\T^2$. In the
non-degenerate case (i.e., if all surfaces $M_\eps$ are in
general position) this curve is of one of three free homotopy
types: $2\times(2,1)$, or $2\times(2,-1)$, or $4\times(1,0)$, see
figure \ref{FigClass}. 

\begin{figure}[ht]
\centerline{\epsfbox{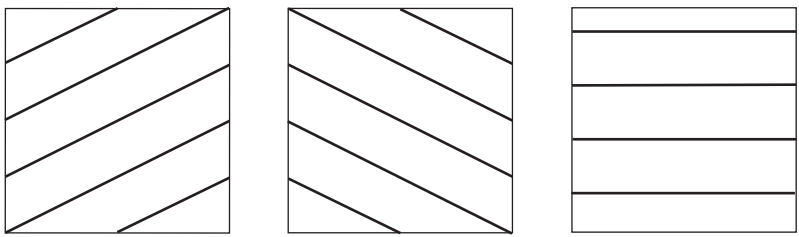}}
\caption{Curves on $\T^2$ of classes 
$2\times(2,1)$, $2\times(2,-1)$ and $4\times(1,0)$.}
\label{FigClass}
\end{figure}

\noindent Indeed, this curve intersects each ``vertical'' cycle
$u=u_0$ in exactly 4 points, while it intersects each
``horizontal'' cycle $v=v_0$ in the same (even) number of points
$\leq2$.

Similarly, the curve $f_{uuu}+f_u=0$ is of one of
three free homotopy types: $2\times(1,2)$, or $2\times(-1,2)$, or
$4\times(0,1)$.

Since the number of intersection points of two
curves of the free homotopy types $n\times(p,q)$ and
$n'\times(p',q')$ is not less than $nn'|pq'-qp'|$, we conclude
that, in our case, this number is at least 12. Indeed, this
number is 12 for the curves $2\times(2,1)$ and $2\times(1,2)$, as
well as for the curves $2\times(2,-1)$ and $2\times(-1,2)$,
 it is equal to 16 in all the cases involving the curves 
$4\times(0,1)$ and $4\times(1,0)$, and it equals 20 for the
intersection of the curves $2\times(-1,2)$ and $2\times(2,-1)$. \proofend
\end{proof}

\begin{remark} 
{\rm 
It is easy to see how the signature of two ``neighbouring''
quadratic points changes, see figure
\begin{figure}[ht]
\centerline{\epsfbox{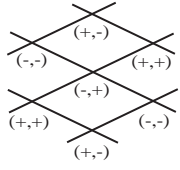}}
\caption{Signature changes.}
\label{PerFig}
\end{figure}
Indeed, the intersecting curves in figure \ref{PerFig} are the curves where $a$
and $b$ change their signs.
The statement then follows from Lemma \ref{ProstLem}.
}
\end{remark}

\medskip

\begin{figure}[ht]
\centerline{\epsfbox{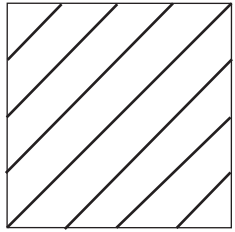}}
\caption{Curves on $\T^2$ of type 
$4\times(1,1)$.}
\label{FigClassBis}
\end{figure}

{\bf B}.
Let us now consider the space of homogeneous second-order
harmonics:
$$ 
f=\cos2u\,(\a_{11}\cos2v+
\a_{12}\sin2v)+
\sin2u\,(\a_{21}\cos2v+
\a_{22}\sin2v),
$$ 
where $\a_{ij}$ are arbitrary constants. In this case, both
curves, $f_{uuu}+f_u=0$, and $f_{vvv}+f_v=0$, on $\T^2$ are
either of homological type $4\times(1,1)$, see figure
\ref{FigClassBis}, or both of type $4\times(-1,-1)$, and they may
avoid intersecting each other altogether. Simple topological
considerations as we used in the proof of Proposition
\ref{OneTwo}, cannot be applied in this case. However, the curves
$f_{uuu}+f_u=0$, and $f_{vvv}+f_v=0$ are no more independent. 

\begin{proposition}
\label{TwoTwo} 
If $f$ is a homogeneous second harmonic then
there are at least 32 distinct points on the
$[0,2\pi)\times[0,2\pi)$-torus at which system (\ref{SSyst}) is
satisfied.
\end{proposition}
\begin{proof} 
It is straightforward to check that the system
(\ref{SSyst}) is equivalent in this case to the following system:
$$
\tau=\frac{\a_{11}\,t+\a_{12}}{\a_{21}\,t+\a_{22}},
\qquad
\tau=\frac{\a_{22}\,t-\a_{21}}{-\a_{12}\,t+\a_{11}},
$$ where $\tau=\tan{2u}$ and $t=\tan{2v}$. This system is
$\frac{\pi}{2}$-periodic and leads to the quadratic equation
$$ (\a_{11}\a_{12}+\a_{21}\a_{22})\,t^2+
(\a_{22}^2-\a_{21}^2+\a_{12}^2-\a_{11}^2)\,t-
(\a_{11}\a_{12}+\a_{21}\a_{22})=0,
$$ whose descriminant is strictly positive. It follows that
system (\ref{SSyst}) has exactly 2 solutions on 
$[0,\frac{\pi}{2})\times[0,\frac{\pi}{2})$ and thus 32 solutions
on $[0,2\pi)\times[0,2\pi)$. \proofend
\end{proof}

\begin{remark} 
{\rm 
Unlike the previous example, the signature of the ``neighbour''
quadratic points is the same, see figure
\ref{NeiBis}. 
\begin{figure}[ht]
\centerline{\epsfbox{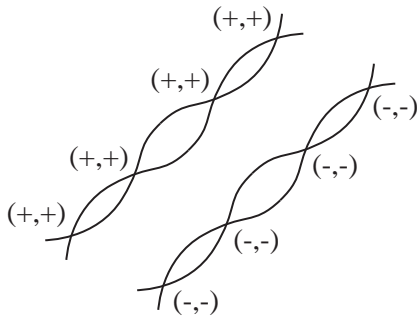}}
\caption{Signature in the homogeneous case.}
\label{NeiBis}
\end{figure}
Indeed, these quadratic points are the points of intersection of the same
curves $a=0$ and $b=0$ which remain transversal to the asimptotic directions,
so that $a$ and $b$ do not change their signs on the corresponding
$x$- and $y$-axes, cf. Lemma \ref{ProstLem}.
}
\end{remark}

\medskip

{\bf C}.
Let us now give an example of a function $f$ for which system
(\ref{SSyst}) has 8 solutions on the
$[0,2\pi)\times[0,2\pi)$-torus. We will consider a sum of
functions of the two above classes. 
\begin{example} {\rm Let $f=\cos(2u-v)+\eps\cos(2u-2v)$. Then the
curve $f_{uuu}+f_u=0$ is of type $2\times(2,1)$, see figure
\ref{FigClass}, while the curve $f_{vvv}+f_v=0$ is of type
$4\times(1,1)$, see figure \ref{FigClassBis}. For $\eps=0$, the
curves intersect transversally, hence, for sufficiently small
$\eps$, the number of intersection points is the same as for
$\eps=0$, that is, equals the number of solutions of the system
$$
\sin (2u-2v)=\sin (2u-v)=0.
$$ This number is equal to 8. }
\end{example}

\section*{Appendix: Wilczynski system of equations} \label{last}

To provide a different description of hyperbolic surfaces in
$\RP^3$ we will write down the system of differential equations
introduced by E. Wilczynski \cite{Wil}\footnote{This reference is
a first systematic study of hyperbolic surfaces in $\RP^3$} (we
also refer to \cite{O-T} for a modern exposition). Given a
parameterized surface $x(u,v)\subset\RP^3$, one wants to lift it
canonically into the vector space $\R^4$ equipped with the
standard volume form.

Let us introduce the notion of
\textit{asymptotic coordinates} in a neighbourhood of an
arbitrary point $m\in{}M$. These are coordinates $(u,v)$ with
origin at $m$ such that the asymptotic lines on $M$ are precisely
the coordinate lines $u=\const$ and $v=\const$. Clearly,
asymptotic coordinates are defined modulo the transformations
$(u,v)\to(U(u),V(v))$. Let first $X(u,v)\subset\R^4$ be an
arbitrary lift. The four vectors $X,X_u,X_v,X_{uv}$ are linearly
independent for every $(u,v)$. One can uniquely fix the lift of
the parameterized surface $x(u,v)$ into $\R^4$ by the condition
\begin{equation}
\label{WilWronCond}
\left| X\,X_u\,X_v\,X_{uv}
\right|=1.
\end{equation} Let us call this lift {\it canonical}.

A straightforward calculation leads to the following fact. The
coordinates of the canonical lift satisfy the system of linear
differential equations
\begin{equation}
\label{WilParform}
\begin{array}{rcl} X_{uu}+a\,X_v+\a\,X & = & 0\\[6pt]
X_{vv}+b\,X_u+\b\,X & = & 0
\end{array}
\end{equation} where $a,b,\a,\b$ are functions in $(u,v)$
satisfying the integrability conditions
\begin{equation}
\label{WilRelEq}
\begin{array}{rcl}
\a_{vv}+b\,\a_u+ 2b_u\,\a-\b_{uu}-2a_v\,\b-a\,\b_v & = & 0\\[4pt]
a\,b_v+2a_v\,b+b_{uu}+2\b_u & = & 0\\[4pt]
b\,a_u+2b_u\,a+a_{vv}+2\a_v & = & 0.
\end{array}
\end{equation} Conversely, system (\ref{WilParform}) whose
coefficients satisfy relations (\ref{WilRelEq}) corresponds to a
non-degenerate parameterized surface
$M\subset\RP^3$.

System (\ref{WilParform}) is called the canonical (or the
Wilczynski) system of differential equations associated with a
surface in $\RP^3$. 
\begin{proposition}
\label{a=b=0Pro} 
The quadratic points on $M$ are the points at
which the functions $a(u,v)$ and $b(u,v)$ vanish.
\end{proposition}
\begin{proof} Identify locally $\RP^3$ and $\R^3$ and consider an
affine lift $X(u,v)$. The linear coordinates $(x,y,z)$ in the
affine 3-space can be chosen in such a way that $X(0,0)$ is the
origin and the vectors $X_u(0,0),X_v(0,0)$ and $X_{uv}(0,0)$ are
the coordinate vectors. Then the surface is locally given by
the equation $z=xy+O(3)$ where $O(3)$ stands for terms, cubic in
$x,y$. One then checks (see, e.g., \cite{O-T}) that the equation
defining
$M$ is
$$
\textstyle z= xy+\frac{1}{3}\,
\left( a\,x^3+b\,y^3
\right) +O(4).
$$ and then applies condition (\ref{CondQuad}). \proofend
\end{proof}

For the sake of completeness let us clarify the geometric meaning of
the coefficients $a$ and $b$.
The proof of the following statement is a straightforward calculation.
\begin{proposition}
\label{VeryLastProp} 
Under coordinate transformations $(u,v)\mapsto(U,V)$, the coefficients
$a$ and $b$ transform as follows:
$$
a(u,v)\mapsto a(U,V)\,\frac{U_u^2}{V_v},
\qquad
b(u,v)\mapsto b(U,V)\,\frac{V_v^2}{U_u}.
$$
\end{proposition}
\noindent
In other words, the following tensor fields
$$
a=a(u,v)\,du^2dv^{-1},
\qquad
b=b(u,v)\,du^{-1}dv^2
$$
are well defined.
Further details can be found in \cite{O-T}, Section 5.1.

\bigskip

{\bf Acknowledgments}. We are pleased to thank C. Duval and V. Kharlamov for their
interest and stimulating discussions. 
The second author was partially supported by an NSF grant.

\end{document}